\documentclass{amsart}
\usepackage{amssymb}
\author[R. Cluckers]{Raf Cluckers$^*$}
\address{Department of Mathematics\\
Katholieke Universiteit Leuven\\ Celestijnenlaan 200B\\ B-3001
Leuven\\ Belgium}
 \email{raf.cluckers@wis.kuleuven.ac.be}
 \urladdr{http://www.wis.kuleuven.ac.be/algebra/raf/}
\thanks{$^*$Research Assistant of the Fund for Scientific Research --
 Flanders (Belgium)(F.W.O.)}
\title{Grothendieck rings of Laurent series fields}
\newtheorem{theorem}{Theorem}
\newtheorem{lemma}{Lemma}
\newtheorem{definition}{Definition}
\newtheorem{proposition}{Proposition}

\newcommand{\Q}{\ensuremath{\mathbb{Q}}}
\newcommand{\C}{\ensuremath{\mathbb{C}}}
\newcommand{\Z}{\ensuremath{\mathbb{Z}}}

\newcommand{\R}{\ensuremath{\mathbb{R}}}
\newcommand{\F}{\ensuremath{\mathbb{F}}}

\newcommand{\Lv}{\ensuremath{\mathcal{L}_{\rm v}}}
\newcommand{\cL}{\ensuremath{\mathcal{L}}}
\newcommand{\Lr}{\ensuremath{\mathcal{L_{\rm ring}}}}
\newcommand{\Lpas}{\cL_{\rm Pas}}
\newcommand{\ac}{\mathrm{ac}}
\renewcommand{\char}{{\rm char}}
\begin{document}
\begin{abstract} We study Grothendieck rings (in the sense of logic)
of fields. We prove the triviality of the Grothendieck rings of
certain fields by constructing definable bijections which imply
the triviality. More precisely, we consider valued fields, for
example, fields of Laurent series  over the real numbers, over
$p$-adic numbers and over finite fields, and construct definable
bijections from the line to the line minus one point.
\end{abstract}
 \maketitle
 \section{introduction}
Recently, the Grothendieck ring of a structure, in the sense of
logic, has been introduced in \cite{DL} and independently in
\cite{KS}. The Grothendieck ring of a model-theoretical structure
is built up as a quotient of the definable sets by definable
bijections (see below), and thus, depends both on the model and
the language. For $(M,\cL)$ a structure with the signature of a
language $\cL$ we write $K_0(M,\cL)$ for the Grothendieck ring of
 $(M,\cL)$. In \cite{CH} and \cite{vdD}, the following explicit
calculations of Grothendieck rings of fields are made:
\begin{quote}
$K_0(\R,\Lr)$ is isomorphic to $\Z$,\\
$K_0(\Q_p,\Lr)$ is trivial,\\
$K_0(\F_p((t)),\Lr)$ is trivial.
\end{quote}
Here, $\Lr$ is the language $(+,-,\cdot,0,1)$. In \cite{DL} and
\cite{KS} it is shown that the Grothendieck ring $K_0(\C,\Lr)$ is
extremely big and complicated; $K_0(\C,\Lr)$, and many other
Grothendieck rings, are not explicitly known.
\par Any Euler characteristic (in the sense of \cite{CH} or
\cite{KS}), defined on the definable sets, factors through the
natural projection of definable sets into the Grothendieck ring,
and, in this sense, to know a Grothendieck ring is to know a
universal Euler characterictic. Nevertheless, it happens that a
Grothendieck ring is trivial.
\par
 The triviality of a
Grothendieck ring can be proven by  constructing a definable
bijection from $X$ to $X\setminus \{a\}$, where $X$ is a definable
set and $\{a\}$ a point on $X$. We develop general techniques to
obtain definable bijections $K\to K^\times$, where $K$ is a valued
field and $K^\times=K\setminus\{0\}$. In section
\ref{sect:laurent} we explain iterated Laurent series fields. In
the present paper we prove:
\begin{theorem}\label{thm:bijection:Qp}
Let $L$ be either $\Q_p$ or a finite field extension of $\Q_p$,
and let $K$ be one of the fields $L$, $L((t_1))$,
$L((t_1))((t_2))$, $L((t_1))((t_2))((t_3))$, and so on. Then
$K_0(K,\Lr)=0$ and there exists a $\Lr$-definable bijection $K\to
K^\times$.\footnote{Here, as always, definable means definable
with parameters.}\footnote{For $K=\Q_p$, this result was proven in
\cite{CH}.}
\end{theorem}
\begin{theorem}\label{thm:bijection:Fq}
Let $L$ be $\F_q$ where $\F_q$ is the finite field with $q=p^l$
elements, $p$ a prime, and let $K$ be one of the fields
$L((t_1))$, $L((t_1))((t_2))$, $L((t_1))((t_2))((t_3))$, and so
on. Then $K_0(K,\Lr)=0$ and there exists a $\Lr$-definable
bijection $K\to K^\times$.\footnote{For $K=\F_q((t))$ this was
proven in \cite{CH}.}
\end{theorem}
Central in the proofs of this paper is a subgroup $H(K,\cL)$ of
$\Z$, associated to a field $K$ and a language $\cL$, which is
sensitive to some elementary arithmetical properties of  the
indices of $n$-th powers in $K^\times$ and of the number of $n$-th
roots in $K^\times$ (see section \ref{sect:calcul}). Using the
definition, it follows immediately that, for example,
\begin{quote}
$H(\R,\Lr)$ is $\Z$,\\
$H(\Q_p,\Lr)$ is $\Z$, and\\
$H(\C,\Lr)$ is $\{0\}$.\\
\end{quote}
We give two criteria for valued fields, for which the value group
has a well-determined minimal strictely positive element, to have
a trivial Grothendieck ring (proposition \ref{prop:HK=Z} and
\ref{prop:crit2}).
\par
We also consider Laurent series fields over $\R$ and over fields
of characteristic $p>0$, using the language of Denef - Pas.
The language of Denef - Pas \cite{Pas} was introduced to study
uniform $p$-adic integrals for all primes $p$, and is now still
used in, for example, the theory of motivic integration (see
\cite{DL} and  \cite{DLinvent}).
\par
For any $\Z$-valued field $K$ with angular component map $\ac$,
the Grothendieck ring $K_0(K,\Lpas)$ is trivial, and there exists
a $\Lpas$-definable bijection from $K^2$ onto $K^2\setminus
\{(0,0)\}$, see \cite{CH}, Thm.~1 and proposition \ref{prop:crit2}
below. (Proposition \ref{prop:crit2} is more general than
\cite{CH}, Thm.~1.) The following theorems give stronger results
for iterated Laurent series fields.
\begin{theorem}\label{thm:bijection:R((t))}
Let $K$ be one of the fields $\R((t_1))$, $\R((t_1))((t_2))$,
$\R((t_1))((t_2))((t_3))$, and so on. We have
\[H(K,\Lr)=\Z.\]
Endow $K$ with a valuation onto a group of the form $\Z^k$ with
lexicographical order, $k>0$, and with the natural angular
component map (as in section \ref{sect:Pas}). Then
\[K_0(K,\Lpas)=\{0\},\]
and there exist a bijection $K\to K^\times$, definable in the
language $\Lpas$ of Denef - Pas.
\end{theorem}
\begin{theorem}\label{thm:bijection3}
Let $L$ be an arbitrary field of characteristic $p>0$. Let $K$ be
one of the fields $L((t_1))$, $L((t_1))((t_2))$,
$L((t_1))((t_2))((t_3))$ and so on. Endow $K$ with a valuation
onto a group of the form $\Z^k$ with lexicographical order, $k>0$,
and with the natural angular component map (as in section
\ref{sect:Pas}). Then $K_0(K,\Lpas)=\{0\}$ and there exists a
bijection $K\to K^\times$ definable in the language of Denef -
Pas.
\end{theorem}
\subsection{Valued fields}
Fix a field $K$. We call $K$ a valued field if there is an ordered
group $(G,+,\leq)$  and a valuation map $v:K\to G\cup\{\infty\}$
such that
\begin{itemize}
\item[(i)] $v(x)=\infty$ if and only if $x=0$;
\item[(ii)] $v(xy)=v(x)+v(y)$ for all $x,y\in K$;
\item[(iii)] $v(x+y)\geq \min\{v(x),v(y)\}$ for all $x,y\in K$.
\end{itemize}
We write $R$ for the  valuation ring $\{x\in K\mid v(x)\geq0\}$ of
$K$, $M$ for its unique maximal ideal and we write $k$ for the
residue field  $R/M$ and $p:R\to k$ for the natural projection. If
$G=\Z$ we call $K$ a $\Z$-valued field. Whenever $K$ is a
$\Z$-valued field, the valuation ring $R$ is a discrete valuation
ring, and a generator $\pi$ of the maximal ideal of $R$ is called
a uniformizer.
\par
A valued field often carries an angular component map modulo $M$,
or angular component map for short; it is a group homomorphism
$\ac:K^\times\to k^\times$, extended by putting $\ac(0)=0$, and
satisfying $\ac(x)=p(x)$ for all $x$ with $v(x)=0$ (see \cite{P}).
\subsection{Iterated Laurent series fields}\label{sect:laurent}
 We define iterated
Laurent series fields by induction. Let $L((t_1))$ be the field of
(formal) Laurent series in the variable $t_1$ over $L$ and let
$L((t_1))\ldots((t_{n-1}))((t_n))$ be the field of (formal)
Laurent series in the variable $t_n$ over
$L((t_1))\ldots((t_{n-1}))$. On a field $L((t_1))\ldots((t_n))$ we
can put many valuations, for example the valuation  $v_n$ taking
values in the lexicographically ordered $n$-fold product of $\Z$,
defined as follows. If $n=1$, then we put $v_1(x)=s\in\Z$ whenever
$x=\sum_{i\geq s}a_st_n^i$ with $a_s\not=0$ and $a_i\in L$. For
general $n$, and $x=\sum_{i\geq s}a_st_n^i$, where $a_s\not=0$ and
$a_i\in L((t_1))\ldots((t_{n-1}))$, we put
$v_n(x)=(s,v_{n-1}(a_s))\in\Z^n$. Remark that the valuation ring
with respect to the valuation $v_n$ is Henselian.
 \subsection{Grothendieck rings}
Let $\cL$ be a language and let $M$ be a model for $\cL$. For
$\cL$-definable sets $X\subset M^m$, $Y\subset M^n$, $m,n>0$, a
$\cL$-definable bijection $X\to Y$ is called an $\cL$-isomorphism
and we write $X\cong_{\cL} Y$, or $X\cong Y$ if the context is
clear, if $X$ and $Y$ are $\cL$-isomorphic. (Definable always
means definable with parameters.) For definable $X$ and $Y$, we
can choose disjoint definable sets $X',Y'\subset K^{m'}$ for some
$m'>0$, such that $X\cong X'$ and $Y\cong Y'$, and then we define
the \emph{disjoint union} $X\sqcup Y$ of $X$ and $Y$ up to
isomorphism as $X'\cup Y'$. By the Grothendieck group $K_0(M,\cL)$
of the structure $(M,\cL)$ we mean the group generated by symbols
$[A]$, for $A$ a $\cL$-definable set, with the relations
$[A]=[A']$ if $A\cong_{\cL} A'$ and $[A]=[B]+[C]$ if $A$ is the
disjoint union of $B$ and $C$. The group $K_0(M,\cL)$ carries a
multiplicative structure induced by $[A\times B]=[A][B]$, where
$A\times B$ is the Cartesian product of definable sets. The
so-obtained ring is called the Grothendieck ring and for a
$\cL$-definable set $X$ we write $[X]$ for the image of $X$ in
$K_0(M,\cL)$.
\par
Let $T$ be a theory in some language $\cL$. A formula $\varphi$
with free variables $x_1,\ldots x_n$ determines a set in $M^n$ for
any model $M$ of $T$. On these sets we can define a disjoint union
operation and Cartesian products in the natural way. The
Grothendieck group $K_0(T,\cL)$ is the group generated by symbols
$[\varphi]$, for $\varphi$ a $\cL$-formula, with the relations
$[\varphi]=[\varphi']$ if the theory $T$ implies that there is
some $\cL$-definable bijection between the sets defined by
$\varphi$ and $\varphi'$, and the relation
$[\varphi]=[\psi]+[\psi']$ if $\varphi$ is the disjoint union of
$\psi$ and $\psi'$.  This group carries a multiplicative structure
induced by the Cartesian product of definable sets and the so
obtained ring is called the Grothendieck ring of $(T,\cL)$.
\section{Languages of Denef - Pas}\label{sect:Pas}
Let $K$ be a valued field, with a valuation map $v:K\to
G\cup\{\infty\}$ for some ordered group $G$, and an angular
component map $\ac:K\to k$, where $k$ is the residue field. Let
$\cL_k$ be an arbitrary expansion of $\Lr$ and let $\cL_G$ be an
arbitrary expansion of the language of ordered groups with
infinity, namely $(+,-,0,\infty,\leq)$. A language of Denef - Pas
can in fact be either language in a wide variety of languages; it
is always a three-sorted language of the form
$(\cL_k,\Lr,\cL_G,v,\ac)$, with as sorts:
\begin{itemize}
\item[(i)] a $k$-sort for the residue field-sort,
\item[(ii)] a $K$-sort for the valued field-sort, and
\item[(iii)] a $G$-sort for the value group-sort.
\end{itemize}
The language $\Lr$ is used for the $K$-sort, $\cL_k$ for the
$k$-sort and $\cL_G$ for the $G$-sort. The function symbol $v$
stands for the valuation map $K\to G\cup\{\infty\}$ and $\ac$
stands for an angular component map $K\to k$ (in fact, this is an
angular component map modulo the maximal ideal $M$). A structure
for a language of Denef - Pas is denoted $(k,K,G\cup\{\infty\})$,
where $k, K$ and $G$ are as above.
\par
Remark that if $G=\Z$, namely if $K$ is a $\Z$-valued field, there
exists a natural angular component map $\ac:K\to k$ sending
$x\not=0$ to $t^{-v(x)}x\bmod M$, where $t$ is a uniformizer of
the valuation ring. More generally, if the value group of $K$ is
$\Z^n$, and $t_1,\ldots,t_n$ are field elements such that
$v(t_1)=(1,0,\ldots,0),\ldots,v(t_n)=(0,\ldots,0,1)$ forms a set
of generators of $\Z^n$, there is a natural angular component map
$\ac:K\to k$ given by $\ac(x)=x\prod_i t_i^{-r_i}\bmod M$, where
$v(x)=(r_1,\ldots,r_n)$. These angular component maps are
canonical up to the choice of $t$ and $t_i$. Languages of Denef -
Pas are denoted $\Lpas$.
\section{Calculations of Grothendieck rings}\label{sect:calcul}
Let $K$ be a field and $\cL$ an expansion of $\Lr$. We  write
$P_n(K)$ or $P_n$ for the $n$-th powers in $K^\times$. For $n>1$
we put
\[r_n(K)=\sharp\{x\in K\mid x^n=1\}\]
and
\[s_n(K)=[K^\times:P_n(K)]
\]
which is either a nonnegative integer or $\infty$.
\begin{definition}\label{def:power} Let $K$ be a field and let
the numbers $s_n(K)$ and $r_n(K)$, $n>1$ be as above. For $n>1$ we
put
\[
\lambda_n(K,\cL)=\frac{s_n(K)}{r_n(K)}
\]
 if the following conditions are satisfied
\begin{itemize}
\item  $s_n(K)<\infty$ and $\frac{s_n(K)}{r_n(K)}\in\Z$ ;\\
\item there exists a $\cL$-definable $n$-th root function. This means that
there exists a definable set $\sqrt[n]{P_n}$ and a definable
bijection $\sqrt[n]{}:P_n(K)\to\sqrt[n]{P_n}$ such that
$(\sqrt[n]{x})^n=x$ for each $x\in P_n(K)$.
\end{itemize}
If one of the above conditions is not satisfied, we put
$\lambda_n(K,\cL)=1$. We define $H(K,\cL)$ as the subgroup of $\Z$
generated by the numbers
\[
\lambda_n(K,\cL)-1
\]
for all $n>1$.
\end{definition}
Remark that if $\cL'$ is an expansion of $\cL$, then there is a
group inclusion $H(K,\cL)\to H(K,\cL')$. Let $\Lv=(\Lr,R)$ be the
language of rings with an extra $1$-ary relation symbol $R$ which
corresponds to a valuation ring inside the model. If the model is
a valued field, we take the natural interpretations.
\begin{lemma}\label{lem:HK}
Let $K$ be a field and $\cL$ an expansion of $\Lr$. For each
positive number $m\in H(K,\cL)$, there exists a $\cL$-definable
bijection
\[
\bigsqcup_{i=1}^{m+1} K^\times\to K^\times,
\]
and thus, in $K_0(K,\cL)$,
\[m [K^\times] =0.
\]
 Moreover, if $K$ is a valued field and $\cL$ is an expansion of $\Lv$, there
exists a $\cL$-definable bijection
\[
\bigsqcup_{i=1}^{m+1} (R\setminus\{0\})\to R\setminus\{0\},
\]
and thus, in $K_0(K,\cL)$,
\[m [R\setminus\{0\}] =0.
\]
\end{lemma}
\begin{proof} We first prove that $K^\times\cong \bigsqcup_{i=1}^{\lambda_n} K^\times$
for all $\lambda_n=\lambda_n(K,\cL)$, $n>1$. If $\lambda_n=1$
there is nothing to prove, so let $\lambda_n>1$. Remark that for
each $x\in K^\times$ and each definable set $A\subset K$ there is
a definable bijection $xA\cong A$. With the notation of definition
\ref{def:power}, the sets $x\sqrt[n]{P_n}$ form a partition of
$K^\times$ when $x$ runs over the $n$-th roots of unity. This
gives $\bigsqcup_{i=1}^{r_n}\sqrt[n]{P_n}\cong K^\times$. Since
$K^\times$ is the disjoint union of all cosets of $P_n$ inside
$K^\times$, we find $\bigsqcup_{i=1}^{s_n} P_n\cong K^\times$.
Combining with the isomorphism $P_n\cong\sqrt[n]{P_n}$ we
calculate:
 \[
 K^\times \cong \bigsqcup_{i=1}^{s_n} P_n \cong \bigsqcup_{i=1}^{s_n} \sqrt[n]{P_n}
\cong
\bigsqcup_{i=1}^{\lambda_n}(\bigsqcup_{i=1}^{r_n}\sqrt[n]{P_n})
  \cong
\bigsqcup_{i=1}^{\lambda_n}K^\times,
  \]
where $s_n=s_n(K)$ and $r_n=r_n(K)$. Now let $m>0$ be in
$H(K,\cL)$ and let $n>1$, $s>0$ be integers. By what we just have
shown, we can add $\lambda_n-1$ disjoint copies of $K^\times$ to
$\bigsqcup_{i=1}^s K^\times$, in the sense that $\bigsqcup_{i=1}^s
K^\times\cong\bigsqcup_{i=1}^{s+\lambda_n-1}K^\times$. Similarly,
if $s>\lambda_n-1$, we can subtract $\lambda_n-1$ disjoint copies
of $K^\times$ from $\bigsqcup_{i=1}^s K^\times$, to be precise,
$\bigsqcup_{i=1}^s
K^\times\cong\bigsqcup_{i=1}^{s-\lambda_n+1}K^\times$. The lemma
follows since the numbers $\lambda_n -1$ generate $H(K,\cL)$.
\par
If $\cL$ is an expansion of $\Lv$, we have the same isomorphisms
and the same arguments for $R\setminus\{0\}$ instead of
$K^\times$, working with $R\cap P_n$ and $R\cap \sqrt[n]{P_n}$
instead of $P_n$ and $\sqrt[n]{P_n}$.
\end{proof}

\begin{definition} Let $R$ be a valuation ring such that the value group
has a minimal strictly positive element. Let $\pi\in R$ have
minimal strictly positive valuation. Write $M$ for the maximal
ideal of $R$. Let $\ac$ be an angular component map $K\to k$,
where $k$ is the residue field. We define the set $R^{(1)}$ as
\[
R^{(1)}=\{x\in R\mid \ac(x)=1\}.
\]
\end{definition}
The set $R^{(1)}$ is not necessarily definable in the language
$\Lr$. If $R^{(1)}$ is definable in some language $\cL$ we have
the following criterion. Remark also that a minimal strictly
positive element in the value group necessarily is unique.
\begin{proposition}\label{prop:HK=Z} Let $K$ be a valued field. Suppose
that the value group has a minimal strictly positive element and
let $\pi\in R$ have this minimal strictly positive valuation. Let
$\cL$ be an expansion of $\cL_{\rm v}$ and let $\ac$ be an angular
component map $K\to k$, with $k$ the residue field. If $R^{(1)}$
is $\cL$-definable and $H(K,\cL)=\Z$, then
\[K_0(K,\cL)=0,\qquad K\cong_\cL K^\times,\ \mbox{ and }\ R\cong_\cL
R\setminus\{0\}.
\]
\end{proposition}
\begin{proof} We first prove that $K_0(K,\cL)=0$. We may suppose
that $\ac(\pi)=1$, otherwise we could replace $\pi$ by $\pi/a$
where $a$ is an arbitrary element with $v(a)=0$ and
$\ac(a)=\ac(\pi)$. The following is a $\cL$-isomorphism
\[
R\sqcup R^{(1)}\to R^{(1)}:\left\{\begin{array}{lcl}x\in R & \mapsto & 1+\pi x,\\
x\in R^{(1)} & \mapsto & \pi x.\end{array}\right.
 \]
This implies, in $K_0(K,\cL)$, that $[R]+[R^{(1)}]=[R^{(1)}]$, and
thus after cancellation, $[R]=0$. By lemma \ref{lem:HK} and
because $H(K,\cL)=\Z$, also $[R\setminus\{0\}]=0$. The following
calculation implies $K_0(K,\cL)=0$:
\[
0=[R]=[R\setminus\{0\}]+[\{0\}]=[\{0\}]=1.
\]
We have $[\{0\}]=1$ because $[\{0\}]$ is the multiplicative unit
in $K_0(K,\cL)$.
 \par
Next we prove $R\cong R\setminus\{0\}$, by taking translates and
applying homotheties to the occurring sets. We make all occurring
disjoint unions explicit. Write $f_1$ for the isomorphism
\[f_1:1+\pi^2(R\setminus\{0\})\to \pi^2(R\setminus\{0\})\cup1+\pi^2(R\setminus\{0\}),\]
given by lemma  \ref{lem:HK}, it is an isomorphism from one copy
of $R\setminus\{0\}$ onto two disjoint copies of
$R\setminus\{0\}$.
Define the function $f_2$ on $\pi^2 R\cup \pi+\pi^2 R^{(1)}$ by
\[f_2:\pi^2 R\cup \pi+\pi^2 R^{(1)}\to \pi+\pi^2 R^{(1)}:
\left\{\begin{array}{rcl}
 \pi^2 x & \mapsto & \pi+\pi^2(1+\pi x),\\
 \pi+\pi^2 x & \mapsto & \pi+\pi^2(\pi x),
\end{array}\right.\]
then $f_2$ is an isomorphism from the disjoint union of $R$ and
$R^{(1)}$ to a copy of $R^{(1)}$. Finally, we find
$\cL$-isomorphisms:
\[f:R\to R\setminus\{0\}:x\mapsto\left\{
\begin{array}{ll}
 f_1(x) & \mbox{if }x\in 1+\pi^2(R\setminus\{0\}),\\
 f_2(x) & \mbox{if }x\in\pi^2 R\cup \pi+\pi^2 R^{(1)},\\
 x      & \mbox{else}
\end{array}\right.\]
and
\[K\to K^\times :x\mapsto\left\{
\begin{array}{ll}
 f(x) & \mbox{if }x\in R,\\
  x   & \mbox{else.}
\end{array}\right.\]
\end{proof}
Proposition \ref{prop:HK=Z} immediately yields the triviality of
the Grothendieck rings of $\Q_p$ and of $\F_q((t))$ with
characteristic different from 2, which was originally proven in
\cite{CH}. Theorems \ref{thm:bijection:Qp} and
\ref{thm:bijection:Fq} of the present paper are generalizations.
In case that $H(K,\cL)$ is different from $\Z$, we formulate the
following criterion. (The argument of this criterion is similar to
the proof of \cite{CH}, Thm.~1.)
\begin{proposition}\label{prop:crit2}
Let $K$ be a valued field. Suppose that the value group has a
unique minimal, strictly positive element and let $\pi\in R$ have
this minimal strictly positive valuation. Let $\cL$ be an
expansion of $\cL_{\rm v}$  and let $\ac$ be an angular component
map $K\to k$. If $R^{(1)}$ is $\cL$-definable, then
\[K_0(K,\cL)=0,\qquad K^2\cong_\cL K^2\setminus\{(0,0)\},\ \mbox{ and }\ R^2\cong_\cL
R^2\setminus\{(0,0)\}.
\]
\end{proposition}
\begin{proof}
We first prove that $K_0(K,\cL)=0$. As above we may suppose that
$\ac(\pi)=1$. The following is a $\cL$-isomorphism
\[
g_1:R\sqcup R^{(1)}\to R^{(1)}:\left\{\begin{array}{lcl}x\in R & \mapsto & 1+\pi x,\\
x\in R^{(1)} & \mapsto & \pi x.\end{array}\right.
 \]
As above, this implies that $[R]=0$ in $K_0(K,\cL)$.
\par
We argument that the disjoint union of two copies of
$(R\setminus\{0\})^2$ is $\cL$-isomorphic to $(R\setminus\{0\})^2$
itself. Define the sets
\begin{eqnarray*}
X_1=\{(x,y)\in (R\setminus\{0\})^2|v(x)\leq v(y)\},\\
X_2=\{(x,y)\in (R\setminus\{0\})^2|v(x)>v(y)\},
\end{eqnarray*}
then $X_1,X_2$ form a partition of $(R\setminus\{0\})^2$. The
isomorphisms
\begin{eqnarray*}
(R\setminus\{0\})^2\to X_1: (x,y)\mapsto (x,xy),\\
(R\setminus\{0\})^2\to X_2: (x,y)\mapsto (\pi xy,y),
\end{eqnarray*}
imply that $(R\setminus\{0\})^2\sqcup (R\setminus\{0\})^2$ is
isomorphic to $X_1\cup X_2$ which is exactly
$(R\setminus\{0\})^2$. After cancellation, it follows that
$[(R\setminus\{0\})^2]=0$.
\par
Since $0=[R]=[R\setminus\{0\}]+[\{0\}]=[R\setminus\{0\}]+1$ we
have $[R\setminus\{0\}]=-1$. Together with
$0=[(R\setminus\{0\})^2]=[R\setminus\{0\}]^2$ this yields $1=0$,
so $K_0(K,\cL)$ is trivial. Write $g_2$ for the isomorphism
\[g_2:(R\setminus\{0\})^2\to (R\setminus\{0\})^2\sqcup (R\setminus\{0\})^2\]
\par
Now take the disjoint union of $R^{(1)}\times (R\setminus\{0\})$,
$R\times (R\setminus\{0\})$ and $(R\setminus\{0\})^2$ inside $R^2$
in some way, meaning that we take disjoint isomorphic copies
inside $R^2$ of the mentioned sets. Using the above isomorphisms
$g_1$ and $g_2$ in a clever way on these disjoint copies, it is
clear that we can remove one copy of $R\times (R\setminus\{0\})$
from $R^2$ and put one copy of $(R\setminus\{0\})^2$ back instead,
hence we find an isomorphism from $R^2$ to itself minus a point.
For details of this construction, we refer to the proof of
\cite{CH}, theorem 1.
\end{proof}
\section{The proofs of theorems 1, 2, and 3}
\begin{proof}[Proof of theorem \ref{thm:bijection:Qp}]
Fix a field $K$ as in the statement. Using Hensel's lemma, it is
elementary to calculate for each $n$ the numbers $r_n(K)$ and
$s_n(K)$, and to find that, for $n$ a prime number,
$s_n(K)/r_n(K)$ is a power of $n$. Further, it is not difficult to
check that taking $n$-th roots is definable. Therefore, the
generator $\lambda_2-1$ of $H(K,\Lr)$ is uneven and $\lambda_3-1$
is even. This implies that $H(K,\Lr)=\Z$.
\par
We calculate explicitly for $K=\Q_p((t))$, for the other fields of
the statement, the arguments are completely similar, although,
notation can get more complicated.
\par
Let $v$ be the valuation on $K$ into $\Z\times\Z\cup\{\infty\}$
with lexicographical order, determined by: $v(x)=(s,r)$ for
$x=\sum_{i\geq s}a_it^i$ with $a_i\in\Q_p$ and $a_s\not=0$, the
$p$-adic valuation of $a_s$ being $r$. The valuation ring $R$ is
definable and can be described  by
\[
R=\{x\in K\mid 1+tx^2\in P_2(K)\ \&\ 1+px^2\in P_2(K)\}
\]
if $p\not=2$ and by
\[
R=\{x\in K\mid 1+tx^3\in P_3(K)\ \&\ 1+px^3\in P_3(K)\}
\]
if $p=2$. Write $M$ for the maximal ideal of $R$. Let $\ac:K\to
\F_p$ be the angular component $\ac(x)=p^{-r}t^{-s}x\bmod M$ for
nonzero $x$ with $v(x)=(s,r)$. The set $R^{(1)}=\{x\in R\mid
\ac(x)=1\}$ is definable since it is the union of the sets
\[
p^it^jP_{p-1}(K),
\]
for $i,j=0,\ldots, p-2$.

\par
Now we can use proposition \ref{prop:HK=Z},
to find a $\Lr$-definable bijection $K\to K^\times$ and to find
that $K_0(K,\Lr)$ is trivial. This proves the proposition for
fields of the form $\Q_p((t))$. When $L$ is a finite field
extension of $\Q_p$ and $K$ an iterated Laurent series field over
$L$, there are $\Lr$-formula's playing the role of $\varphi$ in
the obvious way and the reader can make the adaptations.
\end{proof}
\begin{proof}[Proof of Theorem \ref{thm:bijection:Fq}]
Suppose for simplicity that $K$ is the field
$\big(\F_p((t_1))\big)((t_2))$, where $p$ is a prime. The other
cases are completely similar.
%
Let $v_2$ be the valuation on $K$ into $\Z\times\Z$ as in section
\ref{sect:laurent}; this is a valuation determined by:
$v(x)=(s,r)$ for a Laurent series $\sum_{i\geq s}a_it_2^i$ with
$a_i\in\F_q((t_1))$ and $a_s\not=0$ and $a_s=\sum_{i\geq r}
b_it_1^i$ with $b_r\not=0$ and $b_i\in\F_q$. Write $M$ for the
maximal ideal with respect to $v_2$. The valuation ring $R=\{x\mid
v_2(x)\geq 0\}$ is $\Lr$-definable because of the following
observation:
\[
R=\{x\in K\mid 1+t_2x^2\in P_2(K)\ \&\ 1+t_1x^2\in P_2(K)\},
\]
if $\char(K)\not=2$ and
\[
R=\{x\in K\mid 1+t_2x^3\in P_3\ \&\ 1+t_1x^3\in P_3\},
\]
if $\char (K)=2$.
\par
Let $\ac:K\to \F_p$ be the angular component map $x\mapsto
t_2^{-s}t_1^{-r}x\bmod M$ for nonzero $x$ with $v_2(x)=(s,r)$.
The set $R^{(1)}=\{x\in R\mid \ac(x)=1\}$ is definable since it is
the  union of the sets
\[
t_1^jt_2^iP_{p-1}(K)
\]
for $i,j=0,\ldots,p-2$. Now use proposition \ref{prop:crit2} to
find that $K_0(K,\Lr)$ is trivial and to find a $\Lr$-definable
bijection
\[g_3:R^2\to R^2\setminus\{(0,0)\}.\]
The following is a definable injection:
\[g_4:R^2\to R: (x,y)\mapsto x^p+t_1y^p,\]
and thereby, we can define the $\Lr$-isomorphism
\[
g_5:R\to R\setminus\{0\}: x\mapsto\left\{\begin{array}{ll}
g_4g_3(g_4^{-1}(x)) & \mbox{if } x\in g_4(R^2)\\
x & \mbox{else}.\end{array}\right.
\]
This finishes the proof.
\end{proof}
\begin{proof}[Proof of theorem \ref{thm:bijection:R((t))}]
Let $K$ be the field $\R((t_1))...((t_{n-1}))((t_n))$. Taking
$n$-th roots is clearly $\Lr$-definable, and, using the notation
of definition \ref{def:power} we have that $\lambda_2(K,\Lr)$ is a
power of two and $\lambda_3(K,\Lr)$ is a power of three.
Therefore, $H(\R((t)),\Lr)=\Z$. The existence of a
$\Lpas$-definable bijection $K\to K^\times$ and the triviality of
$K_0(K,\Lpas)$ are formal consequences of proposition
\ref{prop:HK=Z}, because $H(K,\Lr)=\Z$, the value group clearly
has a unique minimal strictly positive element and $R^{(1)}$ is
$\Lpas$-definable.
\end{proof}
\begin{proof}[Proof of theorem \ref{thm:bijection3}]
Let $K$ be $L((t_1))\ldots((t_n))$ where $L$ is a field of
characteristic $p>0$. The statement follows immediately from
proposition \ref{prop:crit2} using the definable injection
\[R^2\to R: (x,y)\mapsto x^p+t_1y^p\]
as in the proof of theorem \ref{thm:bijection:Fq}.
\end{proof}
\bibliographystyle{amsplain}
\bibliography{groth}

\providecommand{\bysame}{\leavevmode\hbox to3em{\hrulefill}\thinspace}
\providecommand{\MR}{\relax\ifhmode\unskip\space\fi MR }
\providecommand{\MRhref}[2]{%
  \href{http://www.ams.org/mathscinet-getitem?mr=#1}{#2}
}
\providecommand{\href}[2]{#2}
\begin{thebibliography}{1}

\bibitem{CH}
Raf Cluckers and Deirdre Haskell, \emph{{G}rothendieck rings of
  $\mathbb{Z}$-valued fields}, Bulletin of Symbolic Logic \textbf{7} (2001),
  no.~2, 262--269.

\bibitem{DLinvent}
J.~Denef and F.~Loeser, \emph{Germs of arcs on singular algebraic varieties and
  motivic integration}, Inventiones Mathematicae \textbf{135} (1999), 201--232.

\bibitem{DL}
\bysame, \emph{Definable sets, motives and $p$-adic integrals}, Journal of the
  American Mathematical Society \textbf{14} (2001), no.~2, 429--469.

\bibitem{KS}
J.~Kraj\'{\i}\v{c}ek and T.~Scanlon, \emph{Combinatorics with definable sets:
  Euler characteristics and grothendieck rings}, Bulletin of Symbolic Logic
  \textbf{6} (2000), 311--330.

\bibitem{Pas}
J.~Pas, \emph{Uniform $p$-adic cell decomposition and local zeta functions},
  Journal f\"ur die reine und angewandte Mathematik \textbf{399} (1989),
  137--172.

\bibitem{P}
\bysame, \emph{On the angular component map modulo $p$}, J. Symbolic Logic
  \textbf{55} (1990), 1125--1129.

\bibitem{vdD}
Lou van~den Dries, \emph{Tame topology and o-minimal structures}, Lecture note
  series, vol. 248, Cambridge University Press, 1998.

\end{thebibliography}
\end{document}